\documentclass[10pt]{amsart}
\usepackage[latin1]{inputenc}
\usepackage{amsmath,latexsym,amssymb,graphicx,dsfont}
\usepackage{natbib}

\newtheorem{Th}{Theorem}[section]
\newtheorem{Lemma}[Th]{Lemma}
\newtheorem{Cor}[Th]{Corollary}
\newtheorem{Prop}[Th]{Proposition}

\newcommand{\E}{\mathbb{E}}
\newcommand{\Prob}{\mathbb{P}}
\newcommand{\calL}{\mathcal{L}}
\newcommand{\calS}{\mathcal{S}}
\newcommand{\calN}{\mathcal{N}}

\numberwithin{equation}{section}

\begin{document}

\title{Rate of Escape on the Lamplighter Tree}
\author[Lorenz A. Gilch]{}

\address{Institut für Mathematik C, University of Technology Graz, Steyrergasse
  30, A-8010 Graz, Austria}

\email{gilch@TUGraz.at}
\urladdr{http://www.math.tugraz.at/$\sim$gilch/}
\date{\today}

\maketitle

\centerline{\scshape Lorenz A. Gilch}
\medskip
{\footnotesize
 \centerline{Graz University of Technology, Graz, Austria}}

\begin{abstract}
Suppose we are given a homogeneous tree $\mathcal{T}_q$ of degree $q\geq 3$,
where at each vertex sits a lamp, which can be switched on or off. This structure can be
described by the wreath product $(\mathbb{Z}/2)\wr \Gamma$,
where $\Gamma=\ast_{i=1}^q \mathbb{Z}/2$ is
the free product group of $q$ factors $\mathbb{Z}/2$. We consider a
transient random walk on a Cayley graph of $(\mathbb{Z}/2)\wr \Gamma$, for which we want to compute lower and upper bounds for the rate of escape, that is, the
speed at which the random walk flees to infinity.
\end{abstract}

\section{Introduction}


Consider a homogeneous tree $\mathcal{T}_q$ of degree $q\geq 3$, where a lamp sits at each vertex, which can have the states 0\,(``off'') or 1\,(``on''). Initially, all
lamps are off. We think of a
lamplighter walking randomly along the tree and switching lamps on or off. Whenever he
stands at a vertex of $\mathcal{T}_q$ he tosses a coin and decides to change the lamp
state at his actual position or to travel to a random neighbour vertex. This is
modeled by a transient
 Markov chain $(Z_n)_{n\geq 0}$, which represents the position of the
 lamplighter and the lamp configuration at time $n$. A natural length function
 $\ell\bigl((\eta,x)\bigr)$, where $\eta$ is a configuration and $x\in\mathcal{T}_q$, is given
 by the length of a shortest path for the lamplighter standing at $x$ to switch all lamps off
 and return to the starting vertex. By transience, our random walk escapes to
 infinity. We are interested in the almost sure, constant limit $\ell=\lim_{n\to\infty}\ell(Z_n)/n$, which describes the
 speed of the random walk. The number $\ell$ is called the
\textit{rate of escape}, or the \textit{drift}. 
It is well-known that the rate of escape exists and is strictly positive for transient random walks on finitely
generated groups. This follows from \textit{Kingman's subadditive
  ergodic  theorem}; see Kingman \cite{kingman}, Derriennic \cite{derrienic} and
Guivarc'h \cite{guivarch}. We provide upper and lower bounds for $\ell$, which
are rather tight. In particular, the random walk escapes faster
to infinity than its projection onto the tree $\mathcal{T}_q$, on which we have the
natural graph metric. In general, the acceleration of the lamplighter random
walk is not obvious. Regarding the case of $\mathcal{T}_2$, Bertacchi \cite{bertacchi} proved that the drift of random
walks on Diestel-Leader graphs and the drift of the random walks' projection
onto $\mathbb{Z}$ coincide. 
\par
Let us briefly review a few selected results regarding the rate of escape.
The classical case is that of random walks on the $d$-dimensional grid $\mathbb{Z}^d$, where $d\geq 1$, which can be described by the sum of $n$
i.i.d. random variables, the increments of $n$ steps. By the law of
large numbers the limit $\lim_{n\to\infty} |Z_n|/n$, where $|\cdot|$ is the
distance on the grid to the starting point of the random walk, exists almost
surely. Furthermore, this limit is positive if the increments have non-zero
mean vector. 
\par
There
are many detailed results for random walks on groups: Lyons,
Pemantle and Peres \cite{lyons-pemantle-peres} gave a lower bound for the rate
of escape of inward-biased random walks on lamplighter groups. Dyubina
\cite{dyubina} proved that the drift on the wreath product $(\mathbb{Z}/2)\wr A$ is
zero, where $A$ is a finitely generated group, if and only if the random walk's projection onto $A$ is recurrent. Revelle \cite{revelle} examined the rate of escape of random walks on wreath
products. He proved laws of the iterated logarithm for the inner and outer
radius of escape. Mairesse
\cite{mairesse-mantaray} computed a explicit formula in terms of the unique
solution of a system of polynomial equations for the rate of escape of
 random walks on the braid group. An important link
between drift and the Liouville property was obtained by Varopoulos
\cite{varopoulos}. He proved that for symmetric finite range random walks on groups the existence of non-trivial bounded
harmonic functions is equivalent to a non-zero rate of escape. This is related
with the link between the rate of escape and the entropy of random walks, compare
e.g. with Kaimanovich and Vershik \cite{kaimanovich-vershik} and Erschler
\cite{erschler2}. The rate of escape has also been studied on trees: Cartwright, Kaimanovich and Woess \cite{cartwright-kaimanovich-woess} investigated the
boundary of homogeneous trees and the drift on them. Nagnibeda and Woess \cite[Section 5]{woess2} proved that the
rate of escape of transient random walks on trees with finitely many cone types is non-zero and give a
formula for it. 
\par
The structure of this article is as follows: In Section \ref{randomwalk} we
explain the structure of the wreath product $(\mathbb{Z}/2)\wr \mathcal{T}_q$,
which encodes our random walk's information, and
define in a natural way a random walk on it. We also sketch the random walk's
convergence behaviour. In Section \ref{low-up} we construct a lower and upper
bound for the rate of escape $\ell$. In Section \ref{low2} we construct another
lower bound for $\ell$, which is in most cases better than the first one. In
Section \ref{remarks} we extend our considerations to two further lamplighter
random walks on trees: Choosing another generating set of $(\mathbb{Z}/2)\wr
\mathcal{T}_q$ and allowing more lamp states, respectively.

\section{Random Walk on the Lamplighter Tree}

\label{randomwalk}
\subsection{The Lamplighter Tree}
Let $3\leq q\in\mathbb{N}$. Consider the homogeneous tree $\mathcal{T}_q$ of degree
$q$, that is, each vertex has $q$ neighbours. Let
$\mathcal{S}:=\{a_1,\dots,a_q\}$. Then all vertices of $\mathcal{T}_q$ can be described
uniquely by finite words over the alphabet $\mathcal{S}$, where no two
consecutive letters are equal, such that we obtain the following symmetric neighbourhood
property: Each $a\in\calS$ is adjacent to the empty word $o$; if
$w\in\mathcal{T}_q$ with last letter $a_i$, then $wa_j$,
$a_j\in\calS\setminus\{a_i\}$, is adjacent to $w$.
 We can
 define a group operation on $\mathcal{T}_q$ by concatenation of words with possible
 cancellations in the middle: if $u,v\in\mathcal{T}_q$ are represented as words over $\calS$, then
 $u\circ v$ is the concatenation with iterated deletions of all blocks of the form
 ``$a_ia_i$''. For instance, if
 $u=a_1a_2a_1$, $v=a_1a_2a_3$, then $u\circ v=a_1a_3$. In particular, the
 identity is $o$ and we have
 $a_i^{-1}=a_i$ for all $i\in\{1,\dots,q\}$. With this defintion
 $\mathcal{T}_q$ is the Cayley graph of the free product group
$\mathbb{Z}/2 \ast \dots \ast\mathbb{Z}/2$ of $q$ factors $\mathbb{Z}/2$, and
in the sequel we shall identify $\mathcal{T}_q$ with this group.
\par
Furthermore, assume that there sits a lamp at each vertex of $\mathcal{T}_q$, which can be switched
off or on, encoded by ``0'' and ``1''. We think of a lamplighter walking along the
tree and switching lamps on and off. The set of finitely supported
configurations of lamps is
$$
\mathcal{N}:=\bigl\lbrace \eta : \mathcal{T}_q\to \mathbb{Z}/2 \ \bigl| \ |\textrm{supp}(\eta)|<\infty \bigr\rbrace.
$$
Denote by $\mathds{1}_o$ the indicator
function on $\mathcal{T}_q$ wrt. $o$ and by $\mathbf{0}$ the
zero function on $\mathcal{T}_q$. Consider now the \textit{wreath product} 
$$
\mathcal{L}_q :=  \Bigl( \sum_{x\in \mathcal{T}_q} \mathbb{Z}/2 \Bigr)  \rtimes \mathcal{T}_q = (\mathbb{Z}/2) \wr \mathcal{T}_q 
$$
of $\mathcal{T}_q$ with the direct sum of copies of $\mathbb{Z}/2$ indexed by
$\mathcal{T}_q$. The elements of $\calL_q$ are pairs of the form $(\eta,x)\in \calN\times
\mathcal{T}_q$, where $\eta$ represents a configuration of the lamps and $x$ the position
of the lamplighter. For \mbox{$x,w\in\mathcal{T}_q$} and $\eta\in\calN$, define
$$
(x\eta)(w):=\eta(x^{-1}w).
$$
The group operation on $\mathcal{L}_q$ is given by
$$
(\eta_1,x)(\eta_2,y) := \bigl( \eta_1 \oplus (x\eta_2),xy\bigr),
$$
where $x,y\in\mathcal{T}_q$, $\eta_1,\eta_2\in \calN$, $\oplus$ is
the componentwise addition modulo 2 and $(\mathbf{0},o)$ is the identity. We
call $\calL_q$ together with this operation the \textsf{Lamplighter Tree}.
\par
Let
$$
\mathcal{S}_{\mathcal{L}_q}:=\bigl\lbrace (\mathds{1}_o,o),(\mathbf{0},a_i)
\mid a_i\in\calS \bigr\rbrace.
$$
Consider the Cayley graph of $\mathcal{L}_q$ with respect to $\mathcal{S}_{\mathcal{L}_q}$. 
We define a length function on $\mathcal{L}_q$ by $\ell\bigl((\eta,x)\bigr)$, which is the length
of the shortest path in the Cayley graph
from $(\eta,x)$ to $(\mathbf{0},o)$. This is the minimal amount of time needed for the
lamplighter to switch off all lamps and walk back to
$o$, when starting at $x$ with configuration $\eta$. Denote by $|x|$ the tree distance of $x\in\mathcal{T}_q$ to $o$ inside $\mathcal{T}_q$.
\par
We now construct a nearest neighbour lamplighter random walk on the wreath
product $\calL_q$. 
Let \mbox{$p\in(0,1)$.} Consider the sequence of i.i.d. random variables
$(\mathbf{i}_k)_{k\in\mathbb{N}}$ valued in $\calL_q$, the increments, with distribution
$$
\mu(w) =
\begin{cases}
p & \textrm{, if } w=(\mathds{1}_o,o)\\
(1-p)/q & \textrm{, if } w=(\mathbf{0},a_i) \textrm{ for some } a_i\in
\mathcal{S}\\
0 & \textrm{, otherwise}
\end{cases}.
$$
A lamplighter random walk starting at $(\mathbf{0},o)$ is described
by $(Z_n)_{n\in\mathbb{N}_0}$ in the
following natural way:
$$
Z_0:=(\mathbf{0},o), \quad Z_n:=Z_{n-1}\mathbf{i}_n \ \textrm{ for all } n\geq 1.
$$
The distribution of $Z_n$ is $\mu^{(n)}$, the $n$-th convolution power of $\mu$
with respect to the group structure of $\mathcal{L}_q$.
More precisely, we write \mbox{$Z_n=(\eta_n,X_n)$}, where $\eta_n$ is the random configuration of the lamps at
time $n$ and $X_n$ is the random vertex at which the lamplighter
stands  at time $n$. We write $\Prob_z[\,\cdot\,]:=\Prob[\,\cdot \mid Z_0=z]$
for any $z\in\calL_q$, if we want to start the lamplighter walk at $z$ instead of
$(\mathbf{0},o)$. We omit this subindex, if we start at $(\mathbf{0},o)$.
\par
Our aim is to estimate the almost sure, constant limit
$$
\ell = \lim_{n\to\infty} \frac{\ell(Z_n)}{n},
$$
which is called \textsf{rate of escape} or \textsf{drift}. Existence of the constant $\ell$
is a consequence of \textit{Kingman's subadditive ergodic theorem}; see Derriennic \cite{derrienic} and
Guivarc'h \cite{guivarch}. It is well-known that simple random walk on
$\mathcal{T}_q$ has rate of escape $(q-2)/q$. Furthermore, we obtain for
our random walk:
\begin{Lemma}
$$
\lim_{n\to\infty} \frac{|X_n|}{n} =  (1-p)\frac{q-2}{q} \quad \Prob-a.s.
$$
\end{Lemma}
\begin{proof}
Standing at some $x\in\mathcal{T}_q\setminus\{o\}$, we move away from $o$ with
probability $(1-p)(q-1)/q$ and towards $o$ with probability
$(1-p)/q$. Thus, $|X_n|$ is a classical birth-and-death Markov chain on
the non-negative integers. Therefore
$$
\lim_{n\to\infty} \frac{|X_n|}{n} = (1-p)\frac{q-1}{q} - \frac{1-p}{q} = (1-p)\frac{q-2}{q}.
$$
\end{proof}
As a consequence, our lamplighter random walk is transient since the projection
$(X_n)_{n\in\mathbb{N}_0}$ onto the tree is transient. 
\par
We now state two lemmata which we will use several times in later
computations. For this purpose, let for $y\in\mathcal{T}_q$ be 
$$
T_y := \min\bigl\lbrace m\geq 1 \mid X_m=y\bigr\rbrace.
$$
the first return stopping time of $y$.
\begin{Lemma}
If $z=(\eta_x,x)\in\calL_q$ and $y\in\mathcal{T}_q$ is a neighbour of $x$ in the tree, then
$$
F:=\Prob_{z}[T_y<\infty] = \frac{1}{q-1}.
$$
\end{Lemma}
\begin{proof}
By vertex-transitivity, it is obvious that $\Prob_{z}[T_y<\infty]$ depends
only on the neighbourhood property and not on the specific points $x$ and $y$. So we get the recursive equation
$$
F = \mu\bigl( (\mathbf{0},a_i)\bigr) +  \mu\bigl(
(\mathds{1}_o,o)\bigr)\cdot F + \sum_{a_j\in \calS\setminus\{a_i\}}
\mu\bigl( (\mathbf{0},a_j)\bigr)\cdot  F^2 
$$
for any  $a_i\in\calS$,
or equivalently,
$$
(1-p)\frac{q-1}{q}\cdot F^2 - (1-p) \cdot F + \frac{1-p}{q} =0.
$$
As $(X_n)_{n\in\mathbb{N}_0}$ is transient, $F<1$ has to be fulfilled. Thus, the
right solution of this quadratic equation is $F = 1/(q-1)$.
\end{proof}
\begin{Lemma}
$$
G := \sum_{n\geq 0} \Prob[X_n=o] =
\frac{q-1}{(1-p)(q-2)} 
$$
\end{Lemma}
\begin{proof}
As
$$
\Prob[T_o<\infty]=\mu\bigl(
(\mathds{1}_o,o)\bigr)+\sum_{a_i\in\calS} \mu\bigl(
(\mathbf{0},a_i)\bigr)\cdot F,
$$
it follows that
$$
G = \sum_{n\geq 0} \Prob[T_o<\infty]^n = \frac{1}{1-
  \Prob[T_o<\infty]} = \frac{q-1}{(1-p)(q-2)}.
$$
\end{proof}

\subsection{Convergence to the Boundary}

Our random walk projects onto the two processes $X_n$ on the tree
$\mathcal{T}_q$ and $\eta_n$ on $\calN$,
of which we can investigate convergence. For $w\in\mathcal{T}_q$ define the
\textit{cone} rooted at $w\in\mathcal{T}_q$ as
$$
C_w:=\bigl\lbrace w'\in\mathcal{T}_q\, \bigl|\, w \textrm{ is prefix of } w'\bigr\rbrace.
$$
The complement $\mathcal{T}_q\setminus C_w$ is denoted by $\overline{C_{w}}$. The set $\partial
\mathcal{T}_q$ consists of all infinite words over $\calS$ with no two equal
consecutive letters and $\partial C_w$ is the subset of $\partial \mathcal{T}_q$
with words starting with prefix $w$. We write $\widehat{C_w}=C_w \cup \partial
C_w$. Then $\widehat{\mathcal{T}_q} = \mathcal{T}_q \cup \partial \mathcal{T}_q$
becomes a compact space, where the topology on $\mathcal{T}_q$ is discrete,
while a
neighbourhood basis of $\tilde w\in\partial \mathcal{T}_q$ is given by all sets $\widehat
C_w$, where $w$ is prefix of $\tilde w$.
\par
A simple and well-known argument shows that $(X_n)_{n\in\mathbb{N}_0}$
converges almost surely to a
random variable $X_\infty$ valued in $\partial\mathcal{T}_q$ in the sense of the
above topology.
\begin{Lemma}\label{equi-distribution}
Let $a\in\calS$. Then
$$
\Prob[X_\infty \textrm{ has first letter } a]=\frac{1}{q}.
$$
\end{Lemma}
\begin{proof}
By conditioning to the last visit in $o$ before finally walking to $a$ with no
consecutive visit to $o$, we obtain
$$
\Prob[X_\infty \textrm{ has first letter } a] = G \cdot \mu\bigl(
(\mathbf{0},a)\bigr)  \cdot
\bigl(1-F\bigr) = \frac{1}{q}.
$$
\end{proof}
Let $\calN^\ast$ be the set of all functions $\eta:\mathcal{T}_q\to
\mathbb{Z}/2$. By transience, each vertex is visited finitely often providing
that the lamp state of each lamp can be flipped finitely often. Thus,
$(\eta_n)_{n\in\mathbb{N}_0}$ converges almost surely pointwise to a random configuration
$\eta_\infty$ valued in $\calN^\ast$.
\par
Later computations require the following probabilities:
\begin{eqnarray*}
\nu_1 &:= & \Prob[a_1\textrm{ is not first letter of } X_\infty , \eta_\infty(C_{a_1})\not\equiv
\mathbf{0}] \quad \textrm{ and}\\
\nu_2 &:= &\Prob[a_1\textrm{ is first letter of } X_\infty,
\eta_\infty(\overline{C_{a_1}}) \equiv
\mathbf{0}].
\end{eqnarray*}
There is a simple relation between $\nu_1$ and $\nu_2$: By vertex-transitivity
and Lemma \ref{equi-distribution}, we have
$$
\nu_1 = F \cdot \Prob[a_1\textrm{ is first letter of } X_\infty,
\eta_\infty(\overline{C_{a_1}})\not\equiv
\mathbf{0}] = \frac{1}{q-1} \cdot \Bigr( \frac{1}{q} - \nu_2\Bigr).
$$
In the next section we will derive a formula for $\ell$ that depends on
$\nu_1$, $\nu_2$ respectively. We will also give lower bounds for these two probabilities
providing upper and lower bounds for $\ell$.

\section{Lower and Upper Bound}
\label{low-up}

In this section we construct a lower and an upper bound for $\ell$. In
particular, we will see that $\ell > \lim_{n\to\infty} |X_n|/n$, that is, the
random walk on $\calL_q$ flees faster to infinity than its projection onto the
tree $\mathcal{T}_q$.
\par
We reformulate our problem for finding a formula for $\ell$. For this
purpose, we apply a technique going back to Furstenberg \cite{furstenberg}, which was used by Ledrappier
\cite[Section 4\,b]{ledrappier} for free groups, and also by the author \cite{gilch1} for free products of groups. 
\par
By Lebesgue's Dominated Convergence Theorem we have
$$
\lim_{n\to\infty} \frac{\mathbb{E}[\ell(Z_n)]}{n} = \ell .
$$
Thus, if we are able to prove convergence of the sequence 
$$
\Bigl( \mathbb{E}[\ell(Z_{n+1})] - \mathbb{E}[\ell(Z_{n})] \Bigr)_{n\in\mathbb{N}}
$$
then its limit must equal $\ell$. We have
$$
\mathbb{E}[\ell(Z_{n})] = \sum_{h\in \calL_q} \ell(h)\, \mu^{(n)}(h) =
\sum_{g,h\in \calL_q} \mu(g)\,
\ell(h)\, \mu^{(n)}(h)
$$
and 
$$
\mathbb{E}[\ell(Z_{n+1})] = \sum_{g,h\in \calL_q} \ell(gh)\,
\mu(g)\, \mu^{(n)}(h).
$$
Thus we obtain
\begin{eqnarray*}
\mathbb{E}[\ell(Z_{n+1})] -\mathbb{E}[\ell(Z_{n})] & = & \sum_{g\in \calL_q} \mu(g) \sum_{h\in \calL_q} \bigl(\ell(gh)-\ell(h)\bigr)\, \mu^{(n)}(h) \\
& = &  \sum_{g\in \calS_{\calL_q}} \mu(g) \int_{\calL_q} \bigl( \ell(gZ_n)
-\ell(Z_n)\bigr)\, d\mathbb{P}.
\end{eqnarray*}
Define the random variables
$$
Y_{g,n} := \ell(gZ_n) -\ell(Z_n)
$$
for any given $g\in\calS_{\calL_q}$ and $n\in\mathbb{N}$. To understand the behaviour of $Y_{g,n}$
for $n\to\infty$, we now investigate differences of the form
$\ell\bigl(g(\eta,x)\bigr)-\ell\bigl((\eta,x)\bigr)$. For this purpose, define for $a\in\calS$
and $\eta\in \calN$ the configurations
$$
\eta_{a}(w) :=
\begin{cases}
\eta (w), & \textrm{if } w \in C_{a}\\
0, & \textrm{otherwise}
\end{cases}
\quad \textrm{ and } \quad
\overline{\eta_{a}}(w) :=
\begin{cases}
\eta (w), & \textrm{if } w \in \overline{C_{a}}\\
0, & \textrm{otherwise}
\end{cases}.
$$
With this notation we have $\eta = \eta_a \oplus \overline{\eta_{a}}$.

\begin{Prop}\label{diff-prop1}
Let $a\in\calS$, $x\in C_a$ and $\eta\in \calN$. Then
$$
\ell\bigl((\mathbf{0},a)(\eta,x)\bigr) - \ell\bigl((\eta,x)\bigr) =
\begin{cases}
1, & \textrm{if } \overline{\eta_{a}} \not\equiv \mathbf{0} \\
-1, & \textrm{if } \overline{\eta_{a}} \equiv \mathbf{0}
\end{cases}.
$$
\end{Prop}
\begin{proof}
Write $x=ay$ with $y\in \overline{C_{a}}$. Since $\eta_a(w)=1$ if and only if
$(a\eta_a)(aw)=1$ for $w\in\mathcal{T}_q$, we obtain
$$
\ell\bigl((\eta,x)\bigr) = \ell\bigl((\overline{\eta_{a}},o)\bigr) +
\ell\bigl((\eta_a,ay)\bigr) = \ell\bigl((\overline{\eta_{a}},o)\bigr) + 1+ 
\ell\bigl((a\eta_a,y)\bigr).
$$
In the last equation we splitted off the necessary walking step from $o$ to $a$
and ``shifted'' $(\eta_a,ay)$ isometrically by multiplying from the left with
$(\mathbf{0},a)$. Observe that $\bigl|(a\eta_a,y)\bigr|$ equals the minimal
distance of a walk starting in $a$, then realizing the configuration $\eta_a$
before finally reaching $ay$. Note also that $a\overline{C_{a}}=C_a$ and  $aC_{a}=\overline{C_{a}}$. See Figure
\ref{Z-to-gZ}.
\begin{figure}[h]
\begin{center}
\includegraphics[width=5cm]{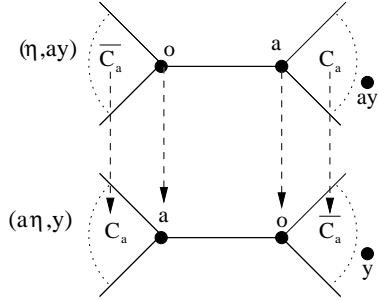}
\end{center}
\caption{Shift from $(\eta,ay)$ to $(a\eta,y)$ with $x\in C_a$}
\label{Z-to-gZ}
\end{figure}
\par
Let $\eta':=a\eta$. Then $(\mathbf{0},a)(\eta,x)=(\eta',y)$. Furthermore,
$\eta'_a=a\overline{\eta_{a}}$ and \mbox{$\overline{\eta'_{a}}=a\eta_a$}. Hence,
$$
\ell\bigl((\eta',y)\bigr) = \ell\bigl((\eta'_{a},o)\bigr) +
\ell\bigl((\overline{\eta'_{a}},y)\bigr) = \ell\bigl((\eta'_{a},o)\bigr) +
\ell\bigl((a\eta_{a},y)\bigr).
$$
As $\overline{\eta_{a}}(w)=1$ if and only if $\eta'_a(aw)=1$, it follows that
$$
\ell\bigl((\eta'_{a},o)\bigr) =
\begin{cases}
2+ \ell\bigl((\overline{\eta_{a}},o)\bigr), & \textrm{if } \overline{\eta_{a}}\not\equiv
\mathbf{0} \\
0, & \textrm{if } \overline{\eta_{a}} \equiv
\mathbf{0}
\end{cases}.
$$
This finishes the proof.
\end{proof}
\begin{Prop}\label{diff-prop2}
Let $a\in\calS$, $x\in \overline{C_{a}}$ and $\eta\in \calN$. Then
$$
\ell\bigl((\mathbf{0},a)(\eta,x)\bigr) - \ell\bigl((\eta,x)\bigr) =
\begin{cases}
-1, & \textrm{if } \eta_{a} \not\equiv \mathbf{0} \\
1, & \textrm{if } \eta_{a} \equiv \mathbf{0}
\end{cases}.
$$
\end{Prop}
\begin{proof}
Observe again that $\eta_a(w)=1$, $\overline{\eta_{a}}(w)=1$ respectively, if and only if
$(a\eta_a)(aw)=1$, $(a\overline{\eta_{a}})(aw)=1$ respectively, for
any $w\in\mathcal{T}_q$.
We obtain
$$
\ell\bigl((\eta,x)\bigr) = \ell\bigl((\eta_{a},o)\bigr) +
\ell\bigl((\overline{\eta_{a}},x)\bigr).
$$
Furthermore,
$$
\ell\bigl((\eta_{a},o)\bigr) =
\begin{cases}
2+ \ell\bigl((a\eta_{a},o)\bigr), & \textrm{if } \eta_a \not\equiv
\mathbf{0}\\
0, & \textrm{if } \eta_a \equiv
\mathbf{0}
\end{cases}.
$$
Let $\eta':=a\eta$. Then $(\mathbf{0},a)(\eta,x)=(\eta',ax)$. Furthermore,
$\eta'_a=a\overline{\eta_{a}}$ and \mbox{$\overline{\eta'_{a}}=a\eta_a$}. See Figure \ref{Z-to-gZ2}.
\begin{figure}[h]
\begin{center}
\includegraphics[width=5cm]{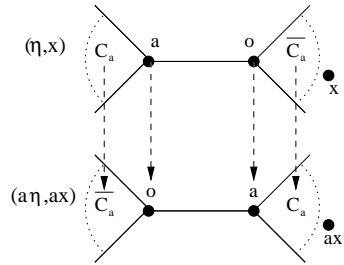}
\end{center}
\caption{Shift from $(\eta,x)$ to $(a\eta,ax)$ with $x\in \overline{C_{a}}$}
\label{Z-to-gZ2}
\end{figure}
\\
Hence,
$$
\ell\bigl((\eta',ax)\bigr) = \ell\bigl((\overline{\eta'_{a}},o)\bigr) +
\ell\bigl((\eta'_{a},ax)\bigr) = \ell\bigl((a\eta_{a},o)\bigr) + 1 +
\ell\bigl((\overline{\eta_{a}},x)\bigr).
$$
This finishes the proof.
\end{proof}
\begin{Prop}\label{diff-prop3}
Let $(\eta,x)\in \calN\times \mathcal{T}_q$. Then
$$
\ell\bigl((\mathds{1}_o,o)(\eta,x)\bigr) - \ell\bigl((\eta,x)\bigr) =
\begin{cases}
1, & \textrm{if } \eta(o) =0 \\
-1, & \textrm{if } \eta(o) =1
\end{cases}.
$$
\end{Prop}
\begin{proof}
Obviously, $(\mathds{1}_o,o)(\eta,x)$ and
$(\eta,x)$ differ only by the lamp state at the root $o$, as $(\mathds{1}_o \oplus \eta)(o)=1-\eta(o)$. This
proves the claim.
\end{proof}
Propositions \ref{diff-prop1},  \ref{diff-prop2} and  \ref{diff-prop3} show
that $Y_{g,n}\in\{-1,1\}$. More precisely, $Y_{g,n}$ remains unchanged after
the last visit in $o$, that is, $Y_{g,n}$ converges almost surely. By Lebesgue's Dominated
Convergence Theorem, almost sure convergence of the sequence
$(\E[\ell(Z_{n+1})-\E[\ell(Z_n)])_{n\in\mathbb{N}}$ follows. Now we want to compute the integrals $\int Y_{g,n}\,d\Prob$. For this purpose, we
need the following probabilities:

\begin{Lemma}\label{prop-onoff}
$$
\Prob[\eta_\infty(o)=0] = \frac{q-2+p}{pq+q-2}\quad \textrm{ and } \quad
 \Prob[\eta_\infty(o)=1] = \frac{p(q-1)}{pq+q-2}.
$$
\end{Lemma}
\begin{proof}
Let 
\begin{eqnarray*}
\widetilde U & := & \Prob[T_o<\infty, X_1\neq o] = \sum_{a\in\calS} \mu\bigl((\mathbf{0},a)\bigr)
F  =  \frac{1-p}{q-1}, \\
\widetilde G & := & \sum_{n\geq 0} \Prob\bigl[X_n=o,\forall j< n: \neg
\bigl(X_j=o \land \mathbf{i}_{j+1}=(\mathds{1}_o,o)\bigr)\bigr] \\
&=& \frac{1}{1-\widetilde
  U} = \frac{q-1}{q-2+p}.
\end{eqnarray*}
Now we can compute the proposed probabilities:
\begin{eqnarray*}
 \Prob[\eta_\infty(o)=0] 
 &=&  \sum_{m\geq 0} \bigl(\widetilde G\cdot p\bigr)^{2m} \cdot \widetilde G
\cdot (1-p) \cdot (1-F)\\
& =& \frac{q-2+p}{pq+q-2},\\
 \Prob[\eta_\infty(o)=1] 
& = & 1- \Prob[\eta_\infty(o)=0]  =  \frac{p(q-1)}{pq+q-2}.
\end{eqnarray*}
\end{proof}
By Propositions \ref{diff-prop1},  \ref{diff-prop2},
\ref{diff-prop3} and Lemma \ref{prop-onoff} we obtain
$$
\int Y_{g,n}\,d\Prob =
\begin{cases}
\frac{(1-p)(q-2)}{pq+q-2}, & \textrm{ if } g=(\mathds{1}_o,o)\\
1-2\nu_1 -2\nu_2, & \textrm{ if } g=(\mathbf{0},a_i) \textrm{ for some } a_i\in\calS
\end{cases}.
$$
Now we can give two explicit formulae for the rate of escape:
\begin{Th}\label{exact-ell}
\begin{eqnarray*}
\ell 
& = & \frac{(1-p)(q-2)}{q} \cdot \Bigl( 1 + 2q \nu_1 +
\frac{pq}{pq+q-2}\Bigr) \\
& = & \frac{(1-p)(q-2)}{q} \cdot \Bigl( \frac{q+1}{q-1} - \frac{2q}{q-1} \nu_2 +
\frac{pq}{pq+q-2}\Bigr)
\end{eqnarray*}
\end{Th}
\begin{proof}
By Lebesgue's Dominated Convergence Theorem and the above computations, we get
\begin{eqnarray*}
\ell &=& \sum_{g\in\calS_{\calL_q}} \mu(g) \int
\lim_{n\to\infty} \bigl( \ell(gZ_{n})-\ell(Z_{n})\bigr)\, d\Prob \\
& = & \sum_{a\in\calS} \Bigl( \mu\bigl((\mathbf{0},a)\bigr) \cdot \bigl(1-
2\nu_1 - 2\nu_2\bigr) \Bigr) +
\mu\bigl((\mathds{1}_o,o)\bigr) \cdot \frac{(1-p)(q-2)}{pq+q-2} \\
&=& (1-p) \cdot (1- 2\nu_1 - 2\nu_2) +  \frac{p(1-p)(q-2)}{pq+q-2}.
\end{eqnarray*}
The rest follows by substituting $\nu_1=\frac{1}{q-1}(\frac{1}{q}-\nu_2)$
resp. $\nu_2=\frac{1}{q}-(q-1)\nu_1$.
\end{proof}
\underline{Remark:}
Observe that $\nu_2=\check G\,
\frac{(1-p)}{q}\, (1-F)$ holds, where
$$
\check G=\sum_{\eta\in\calN'} G(\eta) \quad \textrm{ with} \quad 
G(\eta) = \sum_{n\geq 0} p^{(n)}\bigl((\mathbf{0},o),(\eta,o)\bigr)
$$ 
and $\calN':=\{\eta\in\calN \mid \forall w\in \overline{C_{a_1}}:
\eta(w)=0\}$. The functions $G(\eta)$ are Green functions evaluated at $1$.
 As Green functions are in general hard to compute or even often not
computable and since the structure
of the Cayley graph of $\calL_q$ is very complex, we are only able to give a
lower and upper bound for $\ell$ by estimating $\nu_1$ and $\nu_2$ from
below. For this purpose, we need the following lemma:
\begin{Lemma}
Let $z=(\eta_x,x)\in\calL_q$ and $y\in\mathcal{T}_q$ be a neighbour of $x$ in
the tree. Then the probability that the lamplighter, starting at $x$ with
configuration $\eta_x$,
reaches $y$ without changing any lamps is 
$$
\bar F := \Prob_{z}\bigl[T_y<\infty, \forall k< T_y: \mathbf{i}_k\neq
(\mathds{1}_o,o)\bigr] = \frac{q-\sqrt{q^2-4(q-1)(1-p)^2}}{2(q-1)(1-p)}.
$$
\end{Lemma}
\begin{proof}
By vertex-transitivity, we get the recursive equation
$$
\bar F = \mu\bigl(
(\mathbf{0},a_i)\bigr) + \sum_{a_j\in\calS\setminus\{a_i\}} \mu\bigl(
(\mathbf{0},a_j)\bigr) \bar F^2 \quad \textrm{ for any } a_i\in\calS
$$
with solutions
$$
\bar F = \frac{q \pm \sqrt{q^2-4(q-1)(1-p)^2}}{2(q-1)(1-p)},
$$
where the right one has to to fulfill $\bar F < 1$. This proves the lemma.
\end{proof}
Now we can estimate $\nu_1$ and $\nu_2$ from below:
\begin{Lemma}\label{nu-estimate}
\begin{eqnarray*}
\nu_1  & \geq & \frac{p}{q(pq+q-2)} =: \widehat \nu_1 \quad \textrm{ and}\\
\nu_2  & \geq &
\frac{\widehat G}{1-\widehat G^2p^2}\frac{(1-p)(q-2)}{q(q-1)} =: \widehat\nu_2,
\end{eqnarray*}
where
$$
\widehat G=\frac{2(q-1)}{q-2+\sqrt{q^2-4(q-1)(1-p)^2}}.
$$
\end{Lemma}
\begin{proof}
We restrict the event $[\eta_\infty(C_{a_1})\not\equiv
\mathbf{0}]$ to the event $[\eta_\infty(a_1)=1]$. Thus,
\begin{eqnarray*}
\nu_1 
& \geq & F \cdot \sum_{m\geq 0} \bigl(\widetilde G\cdot p\bigr)^{2m+1} \cdot
\widetilde G \cdot \frac{1-p}{q}\cdot \bigl(1-F\bigr) \\
& = &  \frac{p}{q(pq+q-2)}.
\end{eqnarray*}
For the computation of the lower bound of $\nu_2$, we introduce some further
notation:
\begin{eqnarray*}
\widehat U & := & \Prob\bigl[T_o <\infty, \forall j< T_o:
\neg \bigl( X_j\in \overline{C_{a_1}} \land \mathbf{i}_{j+1} = (\mathds{1}_o,o)\bigr)\bigr]
 \\
& = & \frac{q-1}{q}(1-p) \cdot \bar F + \frac{1-p}{q}\cdot F,\\ 
\widehat G & := & \sum_{n\geq 0} \Prob\bigl[X_n=o, \forall j< n:
\neg \bigl( X_j\in \overline{C_{a_1}} \land \mathbf{i}_{j+1} = (\mathds{1}_o,o)\bigr)\bigr]
\\
&=& \frac{1}{1-\widehat U}.
\end{eqnarray*}
We restrict the event $[\eta_\infty(\overline{C_{a_1}}) \equiv
\mathbf{0}]$ to the event 
that no lamps in $\overline{C_{a_1}}\setminus\{o\}$ are switched on, that is,
$\eta_n(\overline{C_{a_1}}\setminus\{o\})\equiv \mathbf{0}$ for all
$n\in\mathbb{N}$, while we allow to switch the lamp at $o$ for an even number
of switches. This yields
\begin{eqnarray*}
\nu_2
& \geq & \sum_{m\geq 0} \bigl(\widehat G\cdot p\bigr)^{2m} \cdot
 \widehat G \cdot \frac{1-p}{q} \cdot (1-F) \\
& = & \frac{\widehat G}{1-\widehat G^2p^2}\cdot \frac{(1-p)(q-2)}{q(q-1)}.
\end{eqnarray*}
\end{proof}
Now we can give an upper and lower bound for the rate of escape:
\begin{Cor}
\begin{eqnarray*}
\ell & \geq & \frac{(1-p)(q-2)}{q}\cdot \frac{q-2+2p(q+1)}{pq+q-2} =: \ell_{\mathrm{low}} \quad
\textrm{ and}\\
\ell & \leq & \frac{(1-p)(q-2)}{q} \cdot \Bigl( \frac{q+1}{q-1} - \frac{2q}{q-1} \hat \nu_2 +
\frac{pq}{pq+q-2}\Bigr)=: \ell_{\mathrm{up}}
\end{eqnarray*}
\end{Cor}
\begin{flushright}
$\square$
\end{flushright}
Observe that the lower bound also provides $\ell > \lim_{n\to\infty} |X_n|/n$
due to the inequality $(q-2+2p(q+1))/(pq+q-2)>1$, that is, the random walk on $\calL_q$ flees
 to infinity faster than the projection of the random walk onto $\mathcal{T}_q$.
\par
Numerical sample computations are presented at the end of the next section.

\section{Another Lower Bound}

\label{low2}
We construct another lower bound for $\ell$, which is better than $\ell_{\mathrm{low}}$
if $p\leq \frac{q-2}{q-1}$. For this purpose, we give another lower bound for
$\nu_1$, and then apply Theorem \ref{exact-ell}. 
\par
Observe that 
$$
\nu_1 =  F \cdot \underbrace{\Prob[a\textrm{ is first letter of }X_\infty, \eta_\infty(\overline{C_{a_1}})\not\equiv
\mathbf{0}]}_{=:\nu_3}.
$$
Observe that $\eta_\infty(\overline{C_{a_1}})\not\equiv \mathbf{0}$ means that at least one lamp in $\overline{C_{a_1}}$ rests on
forever. Now
we distinguish which of the lamps in $\overline{C_{a_1}}\cap \mathrm{supp}\,\eta_\infty$ is the first
lamp to be switched on and rests finally on, while it is allowed to turn it off
temporarily. More formally, define the random variable $\mathbf{l}_1$ such that
$\mathbf{l}_1=x\in \overline{C_{a_1}}$ if $X_n=X_{n+1}=x$ holds for some $n\in\mathbb{N}$ with
$\eta_m(y)=0$ for all $m<n$ and all $y\in \overline{C_{a_1}}\cap
\mathrm{supp}\,\eta_\infty$. It is sufficient to define $\mathbf{l}_1$ only on
the event $\bigl[\eta_\infty(\overline{C_{a_1}})\not\equiv \mathbf{0}\bigr]$. Define
\begin{eqnarray*}
L & := & \sum_{n\geq 1} \Prob[X_n=a_1,\forall m\in\{1,\dots,n\}: X_m\neq
o] \\
&=& 
\frac{1-p}{q} \cdot \sum_{n\geq 0} \Bigl(\frac{q-1}{q}(1-p) F +p\Bigr)^n
= \frac{1}{q-1}
\end{eqnarray*}
and
$$
\bar G :=  \sum_{n\geq 0} \Prob\bigl[ X_n=o,\forall k\leq n: \mathbf{i}_k\neq
(\mathds{1}_o,o)\bigr] = \frac{1}{1-(1-p)\bar F}.
$$
Now
\begin{eqnarray*}
\nu_3 & = & \sum_{x\in \overline{C_{a_1}}} \Prob[a_1\textrm{ is
  first letter of }X_\infty, \eta_\infty(\overline{C_{a_1}})\not\equiv
\mathbf{0},\mathbf{l}_1=x] \\
& \geq & \sum_{x\in \overline{C_{a_1}}} \bar F^{|x|} \cdot \bar G \cdot
\sum_{m\geq 0} \bigl(p\,\widetilde G\bigr)^{2m+1}\cdot L^{|x|}\cdot
\frac{1-p}{q} \cdot (1-F)\\
& = & \frac{\bar G\widetilde G p}{1-p^2\widetilde G^2}\cdot
\frac{1-p}{q} \cdot \frac{q-2}{q-1} \cdot \sum_{n\geq 0} (q-1)^n \bigl(\bar
F\cdot L\bigr)^n \\
& = & \frac{\bar G \widetilde G p}{1-p^2\widetilde G^2}\cdot
\frac{1-p}{q} \cdot \frac{q-2}{q-1} \cdot \frac{1}{1- \bar F}\\
& = & \frac{p(q-2+p)}{q(pq+q-2)(1-\bar F)(1-(1-p)\bar F)} =:\widehat \nu_3.
\end{eqnarray*} 
Thus,
$$
\ell \geq \frac{(1-p)(q-2)}{q} \cdot \Bigl( 1 + 2\frac{q}{q-1} \widehat\nu_3 +
\frac{pq}{pq+q-2}\Bigr)=\ell_{\mathrm{low},2}.
$$
With the help of \textsc{mathematica} we can show that $\ell_{\mathrm{low},2} \geq
\ell_{\mathrm{low}}$ if $p\leq \frac{q-2}{q-1}$.
\par
Table \ref{compare} compares the values of the trivial lower bound given by 
$\lim_{n\to\infty}|X_n|/n=(1-p)(q-2)/q$, the lower bounds $\ell_{\mathrm{low}}$ and $\ell_{\mathrm{low},2}$ and the upper bound
$\ell_{\mathrm{up}}$ for different values of $q$ and $p$. The \textit{relative
precision of the approximation} is the
quotient 
$$
\frac{\ell_{\mathrm{up}}-\max \{\ell_{\mathrm{low}},\ell_{\mathrm{low},2}\}}{1-\lim_{n\to\infty} \frac{|X_n|}{n}},
$$
which decreases when the degree $q$ of the tree increases: large $q$ yields tighter bounds.
\begin{figure}
\begin{center}
\begin{tabular}[b]{c|c|c|c|c|c|c} 
$q$ & $p$ & $\lim_{n\to\infty} \frac{|X_n|}{n}$ & $\ell_{\mathrm{low}}$ &  $\ell_{\mathrm{low},2}$ & $\ell_{\mathrm{up}}$ & $\substack{\mathrm{relative} \\ \mathrm{precision}}$ \\
&&&&&& \\
\hline
&&&&&& \\
3 & 4/5 & 0.067 & 0.145098 & 0.144410 & 0.157358 & 0.01314\\
3 & 2/3 & 0.111 &  0.234567 & 0.233467 & 0.253778 & 0.02161\\
3 & 1/2 & 0.167 & 0.333 & 0.333 & 0.359733 & 0.03167\\
3 & 1/4 & 0.25  & 0.428571& 0.438050 & 0.461289 & 0.03099 \\
&&&&&& \\
\hline
&&&&&& \\
5 & 4/5	& 0.12 & 0.216 & 0.215942 & 0.221533 &	0.00629\\
5 & 2/3 & 0.2 &	0.347368 & 0.347629 &0.355735 &	0.010459\\
5 & 1/2 & 0.3 &	0.490909 & 0.492585 &0.501825 &	0.01559\\
5 & 1/4 & 0.45 & 0.635294 & 0.641344 & 0.647154 & 0.01056 \\
&&&&&& \\
\hline	
&&&&&& \\
10 & 4/5 & 0.16 & 0.256 & 0.256029 &0.257516 & 0.001805\\
10 & 2/3 & 0.267 & 0.412121 &0.412311 &0.414351 & 0.003040\\
10 & 1/2 & 0.4 &0.584615 &0.585277  &0.587408 &	0.00465 \\
10 & 1/4 & 0.6 & 0.771429 & 0.773099 & 0.774202& 0.00276 \\
&&&&& &\\
\hline			
&&&&&& \\
20 & 4/5 & 0.18 &	0.273176 & 0.273189 &	0.273569 &	0.0004789\\
20& 2/3	 &0.3 &	0.440425 &0.440487 &	0.440994 &	0.0008128\\
20 &	1/2 &	0.45 &	0.626785 &0.626975 &0.627483 &	0.001269\\
20 & 1/4 & 0.675 & 0.836413 & 0.836835 & 0.837079 & 0.00075
\end{tabular}
\end{center}
\caption{Sample computations of lower and upper bounds}
\label{compare}
\end{figure}

\section{Further Random Walk Models}
\label{remarks}

We now consider two other models of lamplighter random walks on $\mathcal{T}_q$
and give lower bounds for the acceleration as compared with their projection
onto the tree.

\subsection{Switch-Walk-Switch}
\label{switch-walk-switch}

Consider again the wreath product $(\mathbb{Z}/2) \wr \mathcal{T}_q$, but now with
generating set
$$
\calS_{\calL_q}^\ast := \bigl\lbrace (\mathds{1}_{A},a) \mid a\in \calS, A \in
\{\emptyset, \{o\},\{a\},\{o,a\}\}\bigr\rbrace.
$$
Consider the random walk on the Cayley graph of $(\mathbb{Z}/2) \wr \mathcal{T}_q$
wrt. $\calS_{\calL_q}^\ast$ described by the sequence of random variables
$\bigl(Z_n\bigr)_{n\in\mathbb{N}_0}$ valued in $(\mathbb{Z}/2)\wr
\mathcal{T}_q$ with $Z_0=(\mathbf{0},o)$, which is governed by the probability measure
$\mu^\ast$ on $\calS_{\calL_q}$ instead of $\mu$, where
$$
\mu^\ast\bigl((\mathds{1}_{A},a_i)\bigr)=
\begin{cases}
\frac{(1-p)^2}{q} & \textrm{, if } A=\emptyset \\
\frac{p(1-p)}{q} &   \textrm{, if } |A|=1 \\
\frac{p^2}{q} &   \textrm{, if } |A|=2 \\
\end{cases}.
$$
This random walk can be interpreted as follows: In one step the lamplighter may
flip the lamp state at his actual position with probability $p$,
walks along one adjacent random edge with probability $1/q$ and may flip the lamp state at the destination
vertex with probability $p$. The number $\ell(Z_n)$ is then the graph distance of $Z_n$ to $(\mathbf{0},o)$ in the
Cayley graph of $(\mathbb{Z}/2)\wr \mathcal{T}_q$
wrt. $\calS^\ast_{\calL_q}$. Write again $Z_n=(\eta_n,X_n)$. Thus,
$(X_n)_{n\in\mathbb{N}_0}$ is simple random walk on $\mathcal{T}_q$.
\par
\begin{figure}[h]
\begin{center}
\includegraphics[width=5cm]{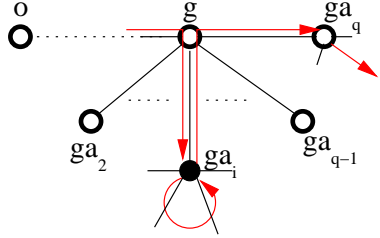}
\end{center}
\caption{Interpretation of $\Delta_k$}
\label{deltak}
\end{figure}

It is well-known that $\lim_{n\to\infty} |X_n|/n=(q-2)/q$. Our aim is to
estimate the ratio of $\ell=\lim_{n\to\infty} \ell(Z_n)/n$ and $(q-2)/q$.
Define for $k\in\mathbb{N}_0$ the exit times
\begin{eqnarray*}
\mathbf{e}_k  & := &  \min \bigl\lbrace m\in\mathbb{N}_0\ \bigl|\ |X_m|=k \land \forall n\geq
m: X_n \in C_{X_m} \bigr\rbrace. 
\end{eqnarray*}
By transience we have almost surely $\mathbf{e}_k<\infty$ for all $k\in\mathbb{N}_0$.
Define now for $k\in\mathbb{N}$ the \textit{pseudo-increments}
$$
\Delta_k:=
\begin{cases}
0, &  \textrm{ if }  \eta_{\mathbf{e}_{k}}(w)=0 \textrm{ for all }
w\in  C_{X_{\mathbf{e}_{k-1}}}\setminus
(C_{X_{\mathbf{e}_{k}}}\cup\{X_{\mathbf{e}_{k-1}}\}) \\
2, &\textrm{ otherwise}
\end{cases}.
$$
The set $C_{X_{\mathbf{e}_{k-1}}}\setminus
(C_{X_{\mathbf{e}_{k}}}\cup\{X_{\mathbf{e}_{k-1}}\})$ is the union of the cones
$C_z$, where $z$ is a forward neighbour of $X_{\mathbf{e}_{k-1}}$ distinct from
$X_{\mathbf{e}_{k}}$. The pseudo-increment $\Delta_k$ represents a lower bound
for the length of a possible deviation inside $C_{X_{\mathbf{e}_{k-1}}}\setminus
C_{X_{\mathbf{e}_{k}}}$, when walking from $o$ to
$X_n$, where $\mathbf{e}_k<n$, with restoring the configuration
$\eta_n$. Note that a
shortest tour from $o$ to $X_n$ does not visit $C_{X_{\mathbf{e}_{k-1}}}\setminus
(C_{X_{\mathbf{e}_{k}}}\cup\{X_{\mathbf{e}_{k-1}}\})$. If at time $\mathbf{e}_{k-1}$ the lamplighter stands at
$g=g'a_1\in\mathcal{T}_q$, then walks to $ga_i$, $i\notin\{1,q\}$, thereby switching the lamp at $ga_i$ on, walks back to $g$ without
flipping the lamp state at $ga_i$, followed by walking to $ga_q$
and rests henceforth in $C_{ga_q}$, then $\Delta_k=2$. See Figure
\ref{deltak}.  
\par
Observe that we have for all $k\geq 1$ 
$$
(*) \quad \quad \ell \bigl((\eta_{\mathbf{e}_k},Z_{\mathbf{e}_k})\bigr) 
\geq  k+\sum_{j=1}^{k}\Delta_j. 
$$
To estimate the distribution of $\Delta_k$, we distinguish if  at
time $\mathbf{e}_{k-1}$ lamps are on in
\mbox{$C_{X_{\mathbf{e}_{k-1}}}\setminus \{X_{\mathbf{e}_k-1}\}$} or not and if lamps are on in \mbox{$C_{X_{\mathbf{e}_{k-1}}}\setminus
\bigl( C_{X_{\mathbf{e}_k}}\cup \{X_{\mathbf{e}_{k-1}}\} \bigr)$} at time $\mathbf{e}_k$.
For $x\in\mathcal{T}_q\setminus\{o\}$ we use the notation $x^-$ to express the
unique neighbour of $x$ \mbox{closer to $o$.} For $k\in\mathbb{N}$ let 
\begin{eqnarray*}
E  & := &   \bigl\lbrace (\eta,x)\in \calN\times (\mathcal{T}_q\setminus\{o\}) \ \bigl| \  \exists w\in C_{x^-}\setminus (C_x\cup \{x^-\}):
\eta(w)=1\bigr\rbrace,\\
E_{k,0} & := & \bigl\lbrace (\eta,x)\in \calN\times \mathcal{T}_q \ \bigl| \ |x|=k,
\forall w\in C_x\setminus \{x\}: \eta(w)=0\bigr\rbrace \quad \textrm{ and}\\
E_{k,2} & := & \bigl\lbrace (\eta,x)\in \calN\times \mathcal{T}_q \ \bigl| \ |x|=k,
\exists w\in C_x\setminus \{x\}: \eta(w)=1\bigr\rbrace.
\end{eqnarray*}
Observe that for
$k\geq 2$ and $r\in\{0,2\}$ it is
$$
\Prob[ Z_{\mathbf{e}_{k-1}}\in E_{k-1,r}] =
\sum_{m\geq 0} \sum_{(\eta,x)\in E_{k-1,r}}
\Prob
\bigl[X_{m-1}=x^-,Z_m=(\eta,x) \bigr]  \cdot
  \bigl(1-F\bigr).
$$
Thus, 
\begin{eqnarray*}
&& \Prob[\Delta_k=2 \mid Z_{\mathbf{e}_{k-1}}\in E_{k-1,r}]\\[1ex]
& = & \frac{1}{\Prob[ Z_{\mathbf{e}_{k-1}}\in E_{k-1,r}]} \sum_{m\geq 0} \sum_{(\eta,x)\in E_{k-1,r}}
\Prob\bigl[X_{m-1}=x^-,Z_m=(\eta,x)
\bigr] \cdot \\
&& \quad \cdot  \Bigl( \sum_{l\geq 1} \Prob_{(\eta,x)}\bigl[\forall \tau\leq
  l:X_\tau\neq x^-,X_{l-1}=x,(\eta_l,X_l)\in E\bigr]\Bigl) \cdot
  \bigl(1-F\bigr)\\
&\geq& \inf_{(\eta,x)\in E_{k-1,r}} \sum_{l\geq 1}
\Prob_{(\eta,x)}\bigl[\forall \tau\leq l:X_\tau\neq x^-,X_{l-1}=x,(\eta_l,X_l)\in E\bigr].
\end{eqnarray*}
Now we can prove:
\begin{Lemma}
\label{delta-lemma}
We have $\mathbb{E}[\Delta_k] \geq  B$ for all $k\in\mathbb{N}$, where 
$$
B:=\frac{4}{q^3}\cdot (q-1)\cdot (q-2)\cdot p\cdot (1-p) >0.
$$
\end{Lemma}
\begin{proof}
Let $k\in\mathbb{N}$. By the above computations we get
$$
\Prob[\Delta_k=2 \mid Z_{\mathbf{e}_{k-1}}\in E_{k-1,0}] \geq 2\cdot (q-1)\cdot
\frac{p(1-p)}{q^2} \cdot \frac{q-2}{q}= \frac{1}{2} B >0
$$
and
$$
\Prob[\Delta_k=2 \mid Z_{\mathbf{e}_{k-1}}\in E_{k-1,2}] \geq \frac{q-2}{q}
\geq \frac{1}{2} B.
$$
Thus, we obtain for $k\geq 2$
\begin{eqnarray*}
\mathbb{E}[\Delta_k] &=&  \Prob[Z_{\mathbf{e}_{k-1}}\in E_{k-1,0}]\cdot
\E[\Delta_k \mid Z_{\mathbf{e}_{k-1}}\in E_{k-1,0}] \\
&& \quad +
\Prob[Z_{\mathbf{e}_{k-1}}\in E_{k-1,2}]\cdot \E[\Delta_k \mid
Z_{\mathbf{e}_{k-1}}\in E_{k-1,2}] \geq  B >0.
\end{eqnarray*}
We have to handle the case $k=1$ separately: here, we have
$\Prob[Z_{\mathbf{e}_{0}}\in E_{0,0}]=1$ and thus 
$$
\E[\Delta_1]  \geq 4\cdot q\cdot \frac{p(1-p)}{q^2} \cdot \frac{q-1}{q} \geq B.
$$
\end{proof}
Now  we want to prove the acceleration on the lamplighter tree:
\begin{Th}
For the \textit{switch-walk-switch} lamplighter random walk, 
$$
\ell \geq \frac{q-2}{q}  \cdot (1+B).
$$
\end{Th}
\begin{proof}
Observe that
$$
\frac{q-2}{q}= \lim_{k\to\infty} \frac{|X_{\mathbf{e}_k}|}{\mathbf{e}_k}
= \lim_{k\to\infty} \frac{|X_{\mathbf{e}_k}|}{k}\frac{k}{\mathbf{e}_k}
=  \lim_{k\to\infty} \frac{k}{\mathbf{e}_k}.
$$
Furthermore,
$$
\ell = \lim_{k\to\infty} \frac{\ell(Z_{\mathbf{e}_k})}{\mathbf{e}_k} = 
 \lim_{k\to\infty} \frac{\ell(Z_{\mathbf{e}_k})}{k}\frac{k}{\mathbf{e}_k} = \frac{q-2}{q}
 \lim_{k\to\infty} \frac{\ell(Z_{\mathbf{e}_k})}{k} \quad
\mathbb{P}-a.s..
$$
As $\ell>0$, the limit $\ell_0=\lim_{k\to\infty} \ell(Z_{\mathbf{e}_k})/k$ exists almost surely and is almost surely constant. We show
now that this limit is greater than $1$. By equation $(*)$
$$
\frac{\ell(Z_{\mathbf{e}_k})}{k} \geq 1 + \frac{1}{k} \sum_{j=1}^k \Delta_j.
$$
Define $D_k:=\frac{1}{k}\sum_{j=1}^k \Delta_j$. Then $0\leq D_k\leq 2$ and
by Lemma \ref{delta-lemma} 
$$
\E[D_k]\geq B>0.
$$
As $\limsup_{k\in\mathbb{N}} D_k = 2 - \liminf_{k\in\mathbb{N}} (2-D_k)$, we
can apply Fatou's Lemma and obtain
\begin{eqnarray*}
\E[\limsup_{k\in\mathbb{N}} D_k] & = & 2 - \int \liminf_{k\in\mathbb{N}} (2-D_k)
\,d\Prob \\
&
\geq & 2 - \liminf_{k\in\mathbb{N}} \int (2-D_k)\, d\Prob = \limsup_{k\in\mathbb{N}} \E[D_k].
\end{eqnarray*}
As $\ell_0\geq 1 + \limsup_{k\in\mathbb{N}} D_k$ we can conclude:
$$
\ell_0 \geq 1 + \E[\limsup_{k\in\mathbb{N}} D_k] \geq 1+
\limsup_{k\in\mathbb{N}} \E[D_k] \geq 1 + B.
$$
This finishes the proof. 
\end{proof}
It is also possible to construct lower and upper bounds for the rate of escape
of this random walk by the technique used in the previous
section. Numerical computations show that those bounds are less tight than in
the case of Section \ref{low-up}, that is, the spread between the bounds is greater.

\subsection{Several Lamp States}

Assume now that there sits a lamp at each vertex of $\mathcal{T}_q$, which can take $r$ different
lamp states including off. These different lamp states are encoded by elements of
$\mathbb{Z}/r$, where $0$ represents the state ``off''. Consider now the wreath product $(\mathbb{Z}/r) \wr
\mathcal{T}_q$ with generating set
$$
\calS^{(r)}_{\calL_q} :=\bigl\lbrace
(k\mathds{1}_o,o),(\mathbf{0},a_i) \mid  k\in\{1,\dots,r-1\}, a_i\in\calS \bigr\rbrace.
$$
Given $p\in(0,1)$. Choose $\alpha_1,\dots,\alpha_{r-1}$ such that
$\sum_{k=1}^{r-1} \alpha_k=p$. Then the corresponding random walk on the
lamplighter tree, where each lamp can take $r$ different lamp states, is the
random walk on the Cayley graph of $(\mathbb{Z}/r) \wr
\mathcal{T}_q$, which is governed by the probability measure $\mu_r$ on $\calS^{(r)}_{\calL_q}$:
$$
\mu_r(z) := 
\begin{cases}
\alpha_k & \textrm{, if } z=(k\mathds{1}_o,o)\\
\frac{1-p}{q} & \textrm{, otherwise}
\end{cases}.
$$
For any $z\in(\mathbb{Z}/r) \wr
\mathcal{T}_q$ it is $\ell(z) = \min\{ n \mid \mu_r^{(n)}(z)>0\}$, where
$\mu_r^{(n)}$ is the $n$-th convolution power of $\mu_r$.
Analogous to Section \ref{switch-walk-switch} we can show that the
corresponding rate of escape $\lim_{n\to\infty} \ell(Z_n)/n$ is strictly
greater than the drift of its projection onto $\mathcal{T}_q$, namely $\lim_{n\to\infty}|X_n|/n=(1-p)(q-2)/q$, where $X_n\in\mathcal{T}_q$ is the random
position of the lamplighter at time $n$.

\bibliographystyle{abbrv}
\bibliography{literatur}

\end{document}